\documentclass[hidelinks]{m2an_mod} 

\usepackage[a-1b]{pdfx}

\usepackage{amsmath}
\usepackage{amssymb}

\newcommand   {\rmref}[1]   {{\rm (\ref{#1})}}

\newcommand   {\fourier}[1] {\widehat{#1}}

\newcommand   {\tnorm}[1]   {|\!|\!|#1|\!|\!|}
		 
\newcommand   {\diff}[1]    {\mathrm{d}#1}

\def \dx      {\,\diff{x}}

\def \dr      {\,\diff{r}}
\def \dt      {\,\diff{t}}

\def \domega  {\,\diff{\omega}}
\def \deta    {\,\diff{\eta}}

\def \e       {\mathrm{e}}
\def \i       {\mathrm{i}}

\def \K       {\mathcal{K}}
\def \F       {\mathcal{F}}

\def \B       {\mathcal{B}}

\def \FH      {\widehat{H}}
\def \X       {X}

\begin{document}

\title{The regularity of electronic wave functions in Barron spaces}
    
\runningtitle{The regularity of electronic wave functions in Barron spaces}

\author{Harry Yserentant}

\address{
Technische Universit\"at Berlin,
Institut f\"ur Mathematik,
10623 Berlin, 
Germany,
\email{yserentant@math.tu-berlin.de}}

\date{} 

\begin{abstract} 
The electronic Schr\"odinger equation describes the motion of 
a finite number of electrons under Coulomb interaction forces 
in the field of a finite number of clamped nuclei. It is 
proved that its solutions for negative eigenvalues, below the 
essential spectrum, lie in the spectral Barron spaces $\B^s$ 
for $s<1$. Examples show that this limit generally cannot be
reached or surpassed.
\end{abstract}

\subjclass{35J10, 35B65, 68T07}

%
%
%
%
%
%
%
%
%
%

\keywords{
Schr\"odinger equation, regularity, Barron spaces
}

\maketitle

\vspace{-4ex}


\section{Introduction}
\label{sec1}

\setcounter{equation}{0}

In these years, the Schr\"odinger equation celebrates its 
centenary. It forms the basis of quantum mechanics and is 
fundamental for the study of atomic and molecular systems. 
Because the nuclei are much heavier than the electrons, 
the electrons follow their motion almost instantaneously. 
Therefore, it is common to separate the motion of the 
nuclei from that of the electrons and to start from the 
electronic Schr\"odinger equation, the equation that 
describes the motion of a finite number of electrons in 
the field of a finite number of clamped nuclei, or in 
other words, to look for the eigenvalues and 
eigenfunctions of the Hamilton operator
\begin{equation}    \label{eq1.1}
H\;=\,{}-\,\frac12\,\sum_{i=1}^N\Delta_i \,-\,
\sum_{i=1}^N\sum_{\nu=1}^K\frac{Z_\nu}{\|x_i-a_\nu\|}\,+\,
\frac12\sum_{\underset{\scriptstyle i\neq j}{i,j=1}}^N 
\frac{1}{\|x_i-x_j\|}
\end{equation}
expressed here in dimensionless form or in atomic units. 
It acts on functions with arguments $x_1,\ldots,x_N$ 
in $\mathbb{R}^3$, which are associated with the positions 
of the considered electrons. The $a_1,\ldots,a_K$ in 
$\mathbb{R}^3$ are the fixed positions of the nuclei and 
the values $Z_\nu>0$ the charges of the nuclei in 
multiples of the electron charge.
The solution space of the equation is the Hilbert space 
$H^1(\mathbb{R}^{3N})$. The multiplication with the
potential 
\begin{equation}    \label{eq1.2}
V(x)\,=\;
-\;\sum_{i=1}^N\sum_{\nu=1}^K\frac{Z_\nu}{\|x_i-a_\nu\|}\,+\,
\frac12\sum_{\underset{\scriptstyle i\neq j}{i,j=1}}^N 
\frac{1}{\|x_i-x_j\|}
\end{equation}
is a bounded mapping from  $H^1(\mathbb{R}^{3N})$ 
to $L_2(\mathbb{R}^{3N})$ as can be shown by means 
of a Hardy inequality. This will play an important 
role later on and allows to express the eigenvalue 
problem in weak form. The weak formulation fixes
the solutions at the singular points of the 
interaction potential (\ref{eq1.2}) and at 
infinity.

The regularity properties of the eigenfunctions 
have been of interest to mathematicians for a 
long time, beginning with Kato's seminal work 
\cite{Kato} almost seventy years ago. Far-reaching 
results on the regularity of the eigenfunctions 
in H\"older spaces can be found in 
\cite{Fournais-Hoffmann-Ostenhof-Sorensen} and 
\cite{Fournais-Hoffmann-Ostenhof-Sorensen_2}.
The study of the regularity in Hilbert spaces 
of mixed derivatives began with Yserentant
\cite{Yserentant_2004}. See Kreusler and 
Yserentant \cite{Kreusler-Yserentant}, 
Meng \cite{Meng} and Yserentant
\cite{Yserentant_2010}, \cite{Yserentant_2011}
for further developments. The regularity 
properties of the solutions in Hilbert spaces 
of mixed derivatives constitute the basis for 
interesting complexity estimates and attempts 
to break the curse of dimensionality 
\cite{Yserentant_2010}.

Inspired by the extensive recent work on the 
approximation properties of neural networks, the 
present paper is dedicated to the regularity of 
the eigenfunctions in Barron spaces. Let 
$B_s(\mathbb{R}^n)$, for any real number $s$, 
be the Banach space that consists of the 
measurable functions 
$u:\mathbb{R}^n\to\mathbb{C}$ with finite norm
\begin{equation}    \label{eq1.3}
\|u\|_{1,s}=\int_{\mathbb{R}^n}(1+\|\omega\|^2)^{s/2}
|u(\omega)|\domega.
\end{equation}
The infinitely differentiable functions with rapidly,
superpolynomially decreasing derivatives of all 
orders, in the following referred to shortly as rapidly 
decreasing functions, form a dense subset of these 
spaces. The spectral Barron space $\B^s(\mathbb{R}^n)$, 
$s\geq 0$, is the space of the functions that possess 
a then also unique representation
\begin{equation}    \label{eq1.4}
u(x)=\left(\frac{1}{\sqrt{2\pi}}\right)^n\!
\int_{\mathbb{R}^n}
\fourier{u}(\omega)\,\e^{\,\i\,\omega\cdot x}\domega
\end{equation}
in terms of a function $\fourier{u}\in B_s(\mathbb{R}^n)$,
their Fourier transform. They are Banach spaces under 
the norm (\ref{eq1.3}) of the Fourier transforms. 
Because every rapidly decreasing function is the Fourier
transform of another rapidly decreasing function, the 
rapidly decreasing functions are dense in 
$\B^s(\mathbb{R}^n)$. The functions in 
$\B^s(\mathbb{R}^n)$ are uniformly continuous and 
vanish at infinity by the Riemann-Lebesgue theorem. 
If $s\geq m$, $m$ a natural number, they are 
$m$-times continuously differentiable and their 
derivatives up to order $m$ vanish at infinity. 

Barron spaces play an increasingly important role in 
the analysis of partial differential equations and 
their numerical solution, especially in high space 
dimensions. An instructive example for this kind of 
research is \cite{Xu}, a work in which a novel 
approach for the numerical solution of elliptic
equations is presented and analyzed. The review 
article \cite{DeRyck-Mishra} gives an overview of 
recent developments from a broader perspective. 
Chen et al. \cite{Chen-Lu-Lu-Zhou} study the 
regularity of the solutions of Schr\"odinger-type 
equations in Barron spaces under rather strong 
conditions on the potential. The work of Dus and 
Ehrlacher \cite{Dus-Ehrlacher} is based on similar 
assumptions. The Coulomb potential (\ref{eq1.2}) 
does not fall into the category of the potentials 
considered in these papers. In fact, one cannot expect 
a high Barron space regularity of the solutions of 
the electronic Schr\"odinger equation, since these 
solutions are in general not differentiable and 
possess typical cusps at the singular points of the 
potential. This can already be seen from the example 
of the (non-normalized) ground state wave function
\begin{equation}    \label{eq1.5}
\psi(x)=\exp(-\|x\|)
\end{equation}
of the hydrogen atom. Its Fourier transform is
\begin{equation}    \label{eq1.6}
\fourier{\psi}(\omega)=
\sqrt{\frac{2}{\pi}}\,\frac{2}{(1+\|\omega\|^2)^2}.
\end{equation}
The norms (\ref{eq1.3}) of $\fourier{\psi}$ are 
finite for $s<1$ and tend like $\sim 1/(1-s)$ 
to infinity as $s$ goes to one. That is, this  
wave function is contained in the spaces 
$\B^s(\mathbb{R}^3)$ for $s<1$, but 
not in $\B^1(\mathbb{R}^3)$.
The message of the present paper is that this behavior 
is no exception. The solutions of the electronic 
Schr\"odinger equation for eigenvalues below 
the ionization threshold $\Sigma\leq 0$ lie in 
the Barron spaces $\B^s(\mathbb{R}^{3N})$, 
$s<1$. The weighted $L_1$-norms 
\begin{equation}    \label{eq1.7}
\|\fourier{\psi}\|_{1,s}=\int_{\mathbb{R}^{3N}}(1+\|\omega\|^2)^{s/2}
|\fourier{\psi}(\omega)|\domega, \quad s<1,
\end{equation}
of their Fourier transforms are finite and 
do not grow faster than $\sim 1/(1-s)$ as 
$s$ goes to one.


\section{Estimates for the Coulomb potential, three space dimensions}
\label{sec2}

\setcounter{equation}{0}

The aim of this and of the next section is to provide 
representations and norm estimates for the Fourier 
transforms of the products of Coulomb potentials with 
in the aforementioned sense rapidly decreasing functions. 
These results are based on an integral representation 
of the Fourier transform of such products. The situation 
is complicated by the fact that the potentials are not 
integrable and have no Fourier transform in the classical 
sense. The starting point of our considerations are the 
following two lemmas, which are borrowed from the theory 
of Riesz potentials and play an important role there;
see Section V.1 in \cite{Stein}, for example. 

\begin{lmm}         \label{lm2.1}
Let $0<\alpha<n$ and $\beta=n-\alpha$. For all rapidly 
decreasing functions $\varphi:\mathbb{R}^n\to\mathbb{C}$
then
\begin{equation}    \label{eq2.1}
\int_{\mathbb{R}^n}\frac{1}{\|x\|^\alpha}\,\varphi(x)\dx
\,=\,\frac{(\sqrt{2\pi})^n}{\gamma(\alpha)}
\int_{\mathbb{R}^n}
\frac{1}{\|\omega\|^\beta}\,\fourier{\varphi}(\omega)\domega
\end{equation}
holds, where the constant $\gamma(\alpha)$ depends 
on the dimension $n$ and is given by
\begin{equation}    \label{eq2.2}
\frac{1}{\gamma(\alpha)}=
\left(\frac{1}{\sqrt{\pi}}\right)^n
\frac{1}{2^\alpha}
\frac{\Gamma(\beta/2)}{\Gamma(\alpha/2)}.
\end{equation}
\end{lmm}

\begin{proof}
The proof starts from the identity
\begin{displaymath}
\int_{\mathbb{R}^n}\frac{1}{\|x\|^\alpha}\,\varphi(x)\dx
\,=\, \frac{1}{\Gamma(\alpha/2)}
\int_0^\infty t^{\alpha/2-1}
\left(\,\int_{\mathbb{R}^n}\varphi(x)\,\e^{-t\,\|x\|^2}\dx\right)\!\dt
\end{displaymath}
that results from the integral representation
\begin{displaymath}
\frac{1}{r^\alpha}=\frac{1}{\Gamma(\alpha/2)}
\int_0^\infty t^{\alpha/2-1}\e^{-t r^2}\dt
\end{displaymath}
by means of the Fubini-Tonelli theorem. If one 
applies the Parseval identity to the inner 
integral and reverses the order of integration 
again, one obtains the representation
\begin{displaymath}
\int_{\mathbb{R}^n}\frac{1}{\|x\|^\alpha}\,\varphi(x)\dx \,=\,
\frac{1}{\Gamma(\alpha/2)}\int_{\mathbb{R}^n}\fourier{\varphi}(\omega)
\left(\,\int_0^\infty t^{\alpha/2-1}\left(\frac{1}{\sqrt{2t}}\right)^n
\exp\left(-\,\frac{\|\omega\|^2}{4t}\,\right)\!\dt\right)\!\domega
\end{displaymath}
of the integral. With the function $h(t)=\|\omega\|^2/(4t)$, 
the inner integrand reads
\begin{displaymath}
-\,\frac{(\sqrt{2})^n}{2^\alpha}\,\frac{1}{\|\omega\|^\beta}\,
h(t)^{\beta/2-1}\e^{-h(t)} h'(t).
\end{displaymath}
The inner integral therefore takes the value
\begin{displaymath}
\frac{(\sqrt{2})^n}{2^\alpha}\,\frac{1}{\|\omega\|^\beta}
\int_0^\infty t^{\beta/2-1}\e^{-t}\dt \,=\,
\frac{(\sqrt{2})^n}{2^\alpha}\,\frac{1}{\|\omega\|^\beta}
\,\Gamma\left(\frac{\beta}{2}\right),
\end{displaymath}
which yields (\ref{eq2.1}) and proves the lemma.
\end{proof}

In the distributional sense, the function
$\fourier{f}(\omega)=c(\alpha)/\|\omega\|^\beta$,
$c(\alpha)=(\sqrt{2\pi})^n/\gamma(\alpha)$, is 
therefore the Fourier transform of the function 
$f(x)=1/\|x\|^\alpha$.
If $u:\mathbb{R}^n\to\mathbb{C}$ is a rapidly 
decreasing function, the functions 
\begin{equation}    \label{eq2.3}
V_\alpha u:x\to\frac{1}{\|x\|^\alpha}\,u(x), 
\quad 0<\alpha<n,
\end{equation}
are integrable. Their Fourier transforms are 
convolution integrals, Riesz potentials.

\begin{lmm}         \label{lm2.2}
If $u:\mathbb{R}^n\to\mathbb{C}$ is a rapidly decreasing 
function, the Fourier transform of the function 
\rmref{eq2.3} is
\begin{equation}    \label{eq2.4}
\fourier{V_\alpha u}(\omega)=
\frac{1}{\gamma(\alpha)}\int_{\mathbb{R}^n}
\frac{1}{\|\omega-\eta\|^\beta}\,\fourier{u}(\eta)\deta,
\end{equation}
where $\gamma(\alpha)$ is the dimension-dependent
constant \rmref{eq2.2}.
\end{lmm}

\begin{proof}
One only needs to insert the function
$\varphi(x)=u(x)\,\e^{-\i\,\omega\cdot x}$
into the equation (\ref{eq2.1}).
\end{proof}

In the following, we restrict ourselves to three 
dimensions and study the product of the Coulomb 
potential
\begin{equation}    \label{eq2.5}
V(x)=\frac{1}{\|x\|}
\end{equation}
with rapidly decreasing functions
$u:\mathbb{R}^3\to\mathbb{C}$. From Lemma~\ref{lm2.2} 
we know the Fourier transform
\begin{equation}    \label{eq2.6}
\fourier{Vu}(\omega)=\frac{1}{2\pi^2}\int_{\mathbb{R}^3}
\frac{1}{\|\omega-\eta\|^2}\,\fourier{u}(\eta)\deta
\end{equation}
of the then integrable functions $Vu$. 
In terms of the integral operator 
\begin{equation}    \label{eq2.7}
\K f(\omega)=\frac{1}{2\pi^2}\int_{\mathbb{R}^3}
\frac{1}{\|\omega-\eta\|^2}\,f(\eta)\deta
\end{equation}
and the operator $\F$ mapping a function $u$ to 
its Fourier transform $\F u=\fourier{u}$,
(\ref{eq2.6}) succinctly reads
\begin{equation}    \label{eq2.8}
\F Vu=\K\F u.
\end{equation}
The key to our further reasoning is the following 
norm estimate for the operator $\K$.

\begin{lmm}         \label{lm2.3}
Let $f:\mathbb{R}^3\to\mathbb{C}$ be a rapidly 
decreasing function, and let $\vartheta$ be  
greater than one. Then
\begin{equation}    \label{eq2.9}
\int_{\mathbb{R}^3}\frac{1}{(1+\|\omega\|^2)^{\vartheta/2}}\,
|\K f(\omega)|\domega
\,\leq\, 
\kappa(\vartheta)\int_{\mathbb{R}^3}|f(\omega)|\domega,
\end{equation}
where the constant $\kappa(\vartheta)$ tends to 
infinity as $\vartheta$ goes to one and is 
given by
\begin{equation}    \label{eq2.10}
\kappa(\vartheta)=
\frac{2\,\Gamma((\vartheta-1)/2)}{\sqrt{\pi}\,\Gamma(\vartheta/2)}.
\end{equation}
\end{lmm}

\begin{proof}
As a rapidly decreasing function, $f$ is integrable
and bounded and $\K f$ therefore bounded. It is
\begin{displaymath}
\int_{\mathbb{R}^3}
\frac{1}{(1+\|\omega\|^2)^{\vartheta/2}}\,
|\K f(\omega)|\domega
\,\leq\, \int_{\mathbb{R}^3}
\frac{1}{(1+\|\omega\|^2)^{\vartheta/2}}\,
\bigg(\frac{1}{2\pi^2}\int_{\mathbb{R}^3}
\frac{1}{\|\omega-\eta\|^2}\,|f(\eta)|\deta
\bigg)\domega.
\end{displaymath}
Tonelli's theorem allows us to reverse 
the order of integration and leads to
\begin{displaymath}
\int_{\mathbb{R}^3}
\frac{1}{(1+\|\omega\|^2)^{\vartheta/2}}\,
|\K f(\omega)|\domega
\,\leq\,
\frac{1}{2\pi^2}\int_{\mathbb{R}^3}\bigg(\,
\int_{\mathbb{R}^3}\frac{1}{\|\omega-\eta\|^2}\,
\frac{1}{(1+\|\omega\|^2)^{\vartheta/2}}\,
\domega\bigg)|f(\eta)|\deta,
\end{displaymath}
where equality holds for real-valued, nonnegative  
$f$. To estimate the inner integral, we subdivide 
the $\mathbb{R}^3$ into the set $\Omega_1$ that 
consists of the $\omega$ for which 
$\|\omega\|\leq\|\omega-\eta\|$ holds and
the set $\Omega_2$ on which
$\|\omega-\eta\|<\|\omega\|$. The integration 
over each of these two subdomains contributes 
a value not greater than
\begin{displaymath}
\int_{\mathbb{R}^3} \frac{1}{\|\omega\|^2}\,
\frac{1}{(1+\|\omega\|^2)^{\vartheta/2}}\,
\domega
\,=\,
4\pi\int_0^\infty\!\frac{1}{(1+r^2)^{\vartheta/2}}\,\dr
\end{displaymath}
to the inner integral.
With the function $h(r)=r^2/(1+r^2)$ 
and the constants $a=1/2$ and $b=(\vartheta-1)/2$,
\begin{displaymath}
\frac{1}{(1+r^2)^{\vartheta/2}} =
\frac{1}{2}\,h(r)^{a-1}(1-h(r))^{b-1}h'(r).
\end{displaymath}
The above integral can therefore be reduced 
to Euler's beta integral
\begin{displaymath}
\int_0^1 t^{a-1}(1-t)^{b-1}\dt
=\frac{\Gamma(a)\Gamma(b)}{\Gamma(a+b)}.
\end{displaymath}
This leads to the constant (\ref{eq2.10}) and 
concludes the proof.
\end{proof}

It is not very difficult to show that $\kappa(\vartheta)$ 
exceeds the best possible constant by a factor of two at 
most, considering functions $f\geq 0$ with supports that 
are concentrated around the origin. In fact, the inner 
integral takes its maximum value at $\eta=0$. A proof 
based on the representation of the integral in spherical 
coordinates can be found in the appendix to this paper.
Therefore, the best possible constant is $\kappa(\vartheta)/2$, 
and $\kappa(-s)/2$ is the norm of $\K$ seen as an operator 
from the space of the rapidly decreasing functions equipped 
with the $L_1$-norm to the space $B_s$, $s<-1$. The constant 
$\kappa(\vartheta)$ possesses the representation
\begin{equation}    \label{eq2.11}
\kappa(\vartheta)=
\frac{2\,q(\vartheta)}{\pi}\,\frac{2}{\vartheta-1},
\quad q(\vartheta)=
\frac{\Gamma(1/2)\,\Gamma((\vartheta+1)/2)}{\Gamma(\vartheta/2)}.
\end{equation}
The function $q$ increases strictly, on the interval 
$1\leq\vartheta\leq 2$ from $q(1)=1$ to $q(2)=\pi/2$.


\section{Estimates for the Coulomb potential, the multi-particle case}
\label{sec3}

\setcounter{equation}{0}

Now we have the tools to study the Fourier transform 
of the product of the interaction potential (\ref{eq1.2}) 
with rapidly decreasing functions and to estimate the 
corresponding norms. Let $V_i$ and $\K_i$, $i=1,\ldots,N$, 
be operators as in the previous section, but now acting 
electron-wise and given by
\begin{equation}    \label{eq3.1}
V_iu(x)=\frac{1}{\|x_i\|}\,u(x), \quad
\K_if(\omega_i,\omega')=\frac{1}{2\pi^2}
\int_{\mathbb{R}^3}
\frac{1}{\|\omega_i-\eta_i\|^2}\,f(\eta_i,\omega')\deta_i,
\end{equation}
where in the definition of the second we have decomposed 
the vectors $\omega\in(\mathbb{R}^3)^N$ into a part 
$\omega_i\in\mathbb{R}^3$ and a residual part
$\omega'\in(\mathbb{R}^3)^{N-1}$. First, we transfer the 
results from the previous section to the $N$-particle case.

\begin{lmm}         \label{lm3.1}
If $f:(\mathbb{R}^3)^N\to\mathbb{C}$ is rapidly 
decreasing, $\K_if$ lies in the spaces 
$B_s(\mathbb{R}^{3N})$, $s<-1$, and 
\begin{equation}    \label{eq3.2}
\|\K_if\|_{1,s}\leq\kappa(-s)\|f\|_{1,0}, 
\end{equation}
where the norms are those from \rmref{eq1.3} and
the constant is given by \rmref{eq2.10}.
\end{lmm}

\begin{proof}
Let $\vartheta=-s$. The functions 
$\omega_i\to f(\omega_i,\omega')$ are rapidly 
decreasing and
\begin{displaymath}    
\int_{\mathbb{R}^3}\frac{1}{(1+\|\omega_i\|^2)^{\vartheta/2}}\,
|\K_i f(\omega_i,\omega')|\domega_i
\,\leq\, 
\kappa(\vartheta)\int_{\mathbb{R}^3}|f(\omega_i,\omega')|\domega_i
\end{displaymath}
by Lemma~\ref{lm2.3}. Since $\|\omega_i\|\leq\|\omega\|$,
this implies the weaker, for our purposes sufficient 
estimate
\begin{displaymath}    
\int_{\mathbb{R}^3}\frac{1}{(1+\|\omega\|^2)^{\vartheta/2}}\,
|\K_i f(\omega_i,\omega')|\domega_i
\,\leq\, 
\kappa(\vartheta)\int_{\mathbb{R}^3}|f(\omega_i,\omega')|\domega_i.
\end{displaymath}
Integration over the remaining variables yields the 
proposition.
\end{proof}

\begin{lmm}         \label{lm3.2}
Let $u:(\mathbb{R}^3)^N\to\mathbb{C}$ be rapidly 
decreasing. The Fourier transform of $V_iu\in L_1$ 
is then
\begin{equation}    \label{eq3.3}
\F V_iu=\K_i\F u.
\end{equation}
\end{lmm}

\begin{proof}
We decompose the full Fourier transform $\F$ into 
the product $\F=\F^i\F'$ of the partial Fourier 
transform
\begin{displaymath}
\F^iu(\omega_i,x')=\left(\frac{1}{\sqrt{2\pi}}\right)^3
\int_{\mathbb{R}^3}u(x_i,x')\,\e^{-\i\omega_i\cdot x_i}\dx_i,
\end{displaymath}
which acts upon the here relevant part $x_i\in\mathbb{R}^3$,
and its correspondingly defined counterpart $\F'$ that acts 
upon the remaining variables. Both map rapidly decreasing 
functions to rapidly decreasing functions. 
The multiplication by $V_i$ and  $\F'$ commute, 
so we have $\F V_i u=\F^i V_i\F' u$ and
$\F V_i u=\K_i\F^i\F' u$ by (\ref{eq2.8}).
\end{proof}

To calculate the Fourier transform of the product of 
a rapidly decreasing function $u$ with the potentials
\begin{equation}    \label{eq3.4}
V_{i,a}(x)=\frac{1}{\|x_i-a\|}, \quad a\in\mathbb{R}^3,
\end{equation}
we introduce the operators $\tau_i(a)$ mapping
a function $f$ to the function $\tau_i(a)f$ 
with the values
\begin{equation}    \label{eq3.5}
\tau_i(a)f(\omega)=\e^{-\i\,\omega_i\cdot a}f(\omega).
\end{equation}
A short calculation using (\ref{eq3.3}) then shows 
that the Fourier transform of $V_{i,a}u$ is
\begin{equation}    \label{eq3.6}
\F V_{i,a}u=\tau_i(a)\K_i\tau_i(-a)\F u.
\end{equation}

Let $Q$ be an orthogonal matrix of corresponding 
dimension and assign to it the operator $Q$ with 
the same name mapping an integrable function $f$ 
to the integrable function $Qf$ with the values
\begin{equation}    \label{eq3.7}
Qf(\omega)=f(Q\omega).
\end{equation}
Since the Fourier transform $\F Qf$ of the function 
$Qf$ is the function $Q\F f$, this operator commutes 
with the Fourier transform. Let $V_{ij}$ be the
potential
\begin{equation}    \label{eq3.8}
V_{ij}(x)=\frac{1}{\|x_i-x_j\|}
\end{equation}
and assign to it the orthogonal matrix $Q_{ij}$
that maps the parts in $x_k\in\mathbb{R}^3$ 
of $x\in(\mathbb{R}^3)^N$ to the vectors
\begin{equation}    \label{eq3.9}
x_i\to\frac{x_i-x_j}{\sqrt{2}},\;\;
x_j\to\frac{x_i+x_j}{\sqrt{2}},\;\;
x_k\to x_k\;\text{for $k\neq i,j$}.
\end{equation}
In operator notation, $V_{ij}=Q_{ij}^{-1}V_jQ_{ij}/\sqrt{2}$. 
With the given commutation relations, the representation
\begin{equation}    \label{eq3.10}
\F V_{ij}u=\frac{1}{\sqrt{2}}\,Q_{ij}^{-1}\K_j Q_{ij}\F u
\end{equation}
of the Fourier transform of the product of the potential 
$V_{ij}$ with a rapidly decreasing function $u$ follows.

Now we are ready, can sum up all the different 
terms, introduce the operator
\begin{equation}    \label{eq3.11}
\K\,=\;
-\;\sum_{i=1}^N\sum_{\nu=1}^K Z_\nu\tau_i(a_\nu)\K_i\tau_i(-a_\nu)
\,+\,
\frac{1}{\sqrt{2}}\sum_{\underset{\scriptstyle i<j}{i,j=1}}^N 
Q_{ij}^{-1}\K_j Q_{ij},
\end{equation}
and summarize our considerations as follows.

\pagebreak

\begin{lmm}         \label{lm3.3}
If $f:\mathbb{R}^{3N}\to\mathbb{C}$ is rapidly 
decreasing, $\K f$ lies in the spaces 
$B_s(\mathbb{R}^{3N})$, $s<-1$, and it is
\begin{equation}    \label{eq3.12}
\|\K f\|_{1,s}\leq C\kappa(-s)\|f\|_{1,0},
\end{equation}
where the constant $C$ depends on the  
total charge of the nuclei and the number of 
the electrons, but not on $s$, and $\kappa(-s)$ 
is defined as in \rmref{eq2.10}. Furthermore,
\begin{equation}    \label{eq3.13}
\F Vu=\K\F u
\end{equation}
for all in the given sense rapidly decreasing 
functions $u:\mathbb{R}^{3N}\to\mathbb{C}$.
\end{lmm}

\begin{proof}
The operators $\tau_i(a)$, $Q_{ij}$, and $Q_{ij}^{-1}$
preserve the norms (\ref{eq1.3}). The estimate thus 
follows from Lemma~\ref{lm3.1} by means of the triangle 
inequality. The commutation relation results from those
for the single parts.
\end{proof}


\section{The regularity of the solutions}
\label{sec4}

\setcounter{equation}{0}

To prove the Barron space regularity, we switch to 
the momentum representation of the Schr\"odinger
equation, which is an equation for the Fourier 
transform of the eigenfunctions. The underlying
Hilbert spaces are no longer the Sobolev spaces 
$H^s$, but the spaces $\FH^s$ of the measurable 
functions $u$ with finite norm given by
\begin{equation}    \label{eq4.1}
\|u\|_{2,s}^2=
\int_{\mathbb{R}^{3N}}(1+\|\omega\|^2)^s|u(\omega)|^2\domega.
\end{equation}
They consist of the $L_2$-Fourier transforms $\F_2u$ 
of the functions $u$ in the Sobolev spaces $H^s$ as 
follows from the Parseval identity. A nonzero 
function $u\in \FH^1$ is an admissible solution of 
the electronic Schr\"odinger equation for the eigenvalue 
$\lambda<0$ in Fourier or momentum representation if 
and only if
\begin{equation}    \label{eq4.2}
u+G(\lambda)\F_2 V\F_2^{-1}u=0.
\end{equation}
The operator $G(\lambda)$ denotes the multiplication 
by the function
\begin{equation}    \label{eq4.3}
G(\lambda)=\frac{2}{\|\omega\|^2-2\lambda}
\end{equation}
and the bounded linear operator $V:H^1\to L_2$ the 
multiplication by the potential (\ref{eq1.2}). The 
aim of this section is to show that the solutions 
of this equation lie in the spaces $B_s$ for $s<1$, 
that is, are Fourier transforms of functions in 
the corresponding Barron spaces $\B^s$. Our theory 
is based on a careful study of the operator
\begin{equation}    \label{eq4.4}
T(\lambda)={}-G(\lambda)\F_2 V\F_2^{-1}
\end{equation}
whose fixed points $u\in\FH^1$, $u\neq 0$, are 
the eigenfunctions for the eigenvalue $\lambda$.

\begin{lmm}         \label{lm4.1}
The operator $T(\lambda)$ maps the functions $u\in\FH^1$ 
to functions in $\FH^2$ that satisfy the estimate
\begin{equation}    \label{eq4.5}
\|T(\lambda)u\|_{2,2}\leq C\,\|u\|_{2,1},
\end{equation}
where the constant $C$ depends on the total charge 
of the nuclei, on the number of electrons, and on  
$\lambda$.
\end{lmm}

\begin{proof}
The estimate is essentially based on the boundedness
of the operator $V:u\to Vu$ from the Sobolev space
$H^1(\mathbb{R}^{3N})$ to $L_2(\mathbb{R}^{3N})$, 
which results from the Hardy inequality
\begin{displaymath}
\int_{\mathbb{R}^3}\frac{1}{\|x\|^2}\,|u(x)|^2\dx
\leq 4\int_{\mathbb{R}^3}\|\nabla u(x)\|^2\dx
\end{displaymath}
for functions $u\in H^1(\mathbb{R}^3)$.
For $u\in \FH^1$, the $L_2$-norm of the function 
$f=\F_2 V\F_2^{-1}u$ thus satisfies an estimate
\begin{displaymath}
\|\F_2 V(\F_2^{-1}u)\|_{2,0}\leq 
c\,\|\F_2(\F_2^{-1}u)\|_{2,1}=c\,\|u\|_{2,1}
\end{displaymath}
since the $L_2$-Fourier transform $\F_2$ preserves 
the $L_2$-norm. The rest of the proof is based on 
the identity
\begin{displaymath}
\|G(\lambda)f\|_{2,2}^2=\;4 \int_{\mathbb{R}^{3N}}
\left(\frac{1+\|\omega\|^2}{\|\omega\|^2-2\lambda}\right)^2
|f(\omega)|^2\domega
\end{displaymath}
and the in the following repeatedly occurring estimate
\begin{displaymath} 
2\;\frac{1+\|\omega\|^2}{\|\omega\|^2-2\lambda}
\leq\max\left(2,-\frac{1}{\lambda}\right).
\end{displaymath}
Inserting $f=\F_2 V\F_2^{-1}u$, the proposition 
follows from the estimate above.
\end{proof}

The most important intermediate step in the proof 
of the regularity theorem is the following estimate, 
first only for rapidly decreasing functions, which 
summarizes our considerations from the last two 
sections.

\begin{lmm}         \label{lm4.2}
If $u:\mathbb{R}^{3N}\to\mathbb{C}$ is a rapidly 
decreasing function, $T(\lambda)u$ lies in the 
$L_1$-based spaces $B_s$, $s<1$, and the 
norms \rmref{eq1.3} of $T(\lambda)u$ can be 
estimated as
\begin{equation}    \label{eq4.6}
\|T(\lambda)u\|_{1,s}\leq C_1\kappa(2-s)\|u\|_{1,0},
\end{equation}
where the constant $C_1$ depends on $\lambda$, 
but is independent of $s$, and $\kappa(2-s)$ 
is given by \rmref{eq2.10}.
\end{lmm}

\begin{proof}
For the functions $f\in B_{s-2}$ one has
\begin{displaymath}
\|G(\lambda)f\|_{1,s}=\,2 \int_{\mathbb{R}^{3N}}
\frac{1+\|\omega\|^2}{\|\omega\|^2-2\lambda}\,
(1+\|\omega\|^2)^{(s-2)/2}\,|f(\omega)|\domega,
\end{displaymath}
from which as in the proof of Lemma~\ref{lm4.1}
\begin{displaymath}
\|G(\lambda)f\|_{1,s} \leq\,
\max\left(2,-\frac{1}{\lambda}\right)\!\|f\|_{1,s-2}
\end{displaymath}
follows. 
The inverse Fourier transform $\F_2^{-1}u$ of a rapidly
decreasing function $u$ is a rapidly decreasing function. 
Since the $L_1$-Fourier transform $\F v$ and the 
$L_2$-Fourier transform $\F_2v$ of functions $v$ in 
the intersection of both spaces coincide, for rapidly 
decreasing functions $u$ we have
\begin{displaymath}
\F_2V\F_2^{-1}u=\F V\F_2^{-1}u=\K\F\F_2^{-1}u=\K u
\end{displaymath}
by Lemma~\ref{lm3.3}. Inserting $f=\F_2V\F_2^{-1}u$ 
into the above estimate, (\ref{eq4.6}) follows 
from (\ref{eq3.12}).
\end{proof}

To proceed, we need to introduce a new scale of 
functions spaces $\X_s$, the intersections of 
the space $\FH^1$ and the spaces $B_s$ from the 
introduction. We equip them with the norms
\begin{equation}    \label{eq4.7}
\tnorm{u}_s=\|u\|_{2,1}+\|u\|_{1,s}
\end{equation}
that are composed of the norm on the solution 
space $\FH^1$ of the equation and the weighted 
$L_1$-norms (\ref{eq1.3}).

\pagebreak

\begin{lmm}         \label{lm4.3}         
The spaces~$\X_s$ are complete and the rapidly 
decreasing functions lie dense in them. 
\end{lmm}

\begin{proof}
The proof of the completeness of the spaces is almost 
identical to the proof of the completeness of $L_1$ 
itself. Let $u_1,u_2,\ldots$ be a Cauchy sequence of 
functions in $\X_s$ and let $v_1,v_2,\ldots$ be a 
subsequence for which
\begin{displaymath}
\tnorm{v_{k+1}-v_k}_s\leq 2^{-k}
\end{displaymath}
holds. The monotonously increasing sequence of the 
nonnegative measurable functions $f_1,f_2,\ldots$ 
given by
\begin{displaymath}
f_n(x)=\sum_{k=1}^n|v_{k+1}(x)-v_k(x)| 
\end{displaymath}
then tends pointwise to a measurable limit function
$f$. In the next step the structure of the norms
enters. Because the norm of a function in $X_s$ 
coincides with the norm of its absolute value, 
the norms 
\begin{displaymath}
\tnorm{f_n}_s\leq\sum_{k=1}^n\tnorm{v_{k+1}-v_k}_s
\leq\sum_{k=1}^n 2^{-k} 
\end{displaymath}
remain uniformly bounded and by the monotone 
convergence theorem the norm 
\begin{displaymath}
\tnorm{f}_s=\lim_{n\to\infty}\tnorm{f_n}_s
\end{displaymath}
of the limit function $f$ is finite. This means 
that $f$ is almost everywhere finite, that the 
$v_k$ converge almost everywhere pointwise to 
a measurable function $v$, and that the absolute 
values of the $v_k$ and of $v$ are almost everywhere 
bounded by the function $|v_1|+f$ of finite norm. 
So the dominated convergence theorem leads to
\begin{displaymath}
\lim_{k\to\infty}\tnorm{v-v_k}_s=0.
\end{displaymath}
That is, $v$ is an accumulation point of the 
Cauchy sequence of the $u_k$. The $u_k$ thus 
converge to $v$. 

The functions in $\X_s$ can be approximated by 
functions with bounded support, as again follows
from the dominated convergence theorem, and 
these by infinitely differentiable functions 
with bounded support via the convolution with 
appropriate mollifiers. This shows that the 
rapidly decreasing functions are dense 
in $\X_s$.
\end{proof}

The estimate from Lemma~\ref{lm4.2} can 
be transferred to the functions in 
$\X_0=\FH^1\cap L_1$.

\begin{lmm}         \label{lm4.4}
If $u\in\X_0$, the function $T(\lambda)u$ lies in 
$\X_s$ for $s<1$ and its norms \rmref{eq1.3} 
can be estimated as
\begin{equation}    \label{eq4.8}
\|T(\lambda)u\|_{1,s}\leq C_1\kappa(2-s)\|u\|_{1,0},
\end{equation}
with the same constant $C_1$ as in \rmref{eq4.6}
that is independent of $s$ and $\kappa(2-s)$ as 
in \rmref{eq2.10}.
\end{lmm}

\begin{proof}
Let $u_1,u_2,\ldots$ be a sequence of rapidly 
decreasing functions that converges in $\X_0$ to 
$u$. By Lemma~\ref{lm4.1} and Lemma~\ref{lm4.2}, 
the functions $T(\lambda)u_k$ form Cauchy 
sequences in the given spaces $\X_s$. Since the 
spaces $\X_s$ are complete, the $T(\lambda)u_k$ 
converge in $\X_s$ to a function $f_s$. Since 
convergence in $\X_s$ implies convergence 
in $B_s$,
\begin{displaymath}  
\|f_s\|_{1,s}\leq C_1\kappa(2-s)\|u\|_{1,0}
\end{displaymath}
follows from (\ref{eq4.6}). But convergence in 
$\X_s$ also implies convergence in $\FH^1$. 
By Lemma~\ref{lm4.1}, the $T(\lambda)u_k$
converge in $\FH^1$ and even in $\FH^2$ to 
$T(\lambda)u$. This means 
$\|f_s-T(\lambda)u\|_{2,1}=0$ and 
proves the proposition.
\end{proof}

In the following we split the considered functions 
into a low-frequency and a high-frequency part. 
Let $\Omega>0$ be a given bound that separates 
the low from the high frequencies and let 
$\chi(\omega)=0$ for $\|\omega\|<\Omega$ and 
$\chi(\omega)=1$ for $\|\omega\|\geq\Omega$. Let 
the operator $P$ denote the multiplication with 
the cut-off function $\chi$. The high-frequency 
part $v=Pu$ of the solution $u\in \FH^1$ of the 
equation (\ref{eq4.2}) satisfies the equation
\begin{equation}    \label{eq4.9}
v-PT(\lambda)v=PT(\lambda)(u-Pu).
\end{equation}
This equation is the key to our regularity result, 
which will follow from the contractivity of 
$PT(\lambda)$ seen as an operator from $\FH^1 $ 
to $\FH^1$ as well as from $\X_0$ to $\X_0$ for 
sufficiently large chosen frequency bounds $\Omega$.

\begin{lmm}         \label{lm4.5}
For $u\in\FH^1$, the function $PT(\lambda)u\in\FH^2$
can be estimated as
\begin{equation}    \label{eq4.10}
\|PT(\lambda)u\|_{2,1}\leq 
C\left(\frac{1}{1+\Omega^2}\right)^{1/2}\!\|u\|_{2,1},
\end{equation}
where the constant $C$ is the same as in 
\rmref{eq4.5}.
\end{lmm}

\begin{proof}
This follows from the almost trivial estimate
\begin{displaymath}
\|Pf\|_{2,1}^2\leq\frac{1}{1+\Omega^2}\,\|f\|_{2,2}^2
\end{displaymath}
for functions $f\in\FH^2$ and Lemma~\ref{lm4.1}
applied to $f=T(\lambda)u$.
\end{proof}

\begin{lmm}         \label{lm4.6}
If $u\in\X_0$, the function $PT(\lambda)u$ lies in 
$\X_s$ for $s<1$ and its norms \rmref{eq1.3} 
can be estimated as
\begin{equation}    \label{eq4.11}
\|PT(\lambda)u\|_{1,s}\leq 
C_1\left(\frac{1}{1+\Omega^2}\right)^{\vartheta/2}
\!\kappa(2-(s+\vartheta))\,\|u\|_{1,0},
\end{equation}
where $\vartheta<1-s$ can be chosen arbitrarily and
$C_1$ is the same constant as in \rmref{eq4.6} and
\rmref{eq4.8}.
\end{lmm}

\begin{proof}
The proposition follows from the estimate
\begin{displaymath}
\|Pf\|_{1,s}\leq
\left(\frac{1}{1+\Omega^2}\right)^{\vartheta/2}\!\|f\|_{1,s+\vartheta}
\end{displaymath}
for functions $f\in B_{s+\vartheta}$ and 
from Lemma~\ref{lm4.4} applied to $
f=T(\lambda)u$.
\end{proof}

\begin{thrm}
The solutions $u\in\FH^1$ of the electronic Schr\"odinger 
equation \rmref{eq4.2} in momentum representation for 
eigenvalues $\lambda<0$ lie in the spaces $\X_s$ for
$0\leq s<1$. For these $s$, their norms \rmref{eq1.3}
satisfy the estimate
\begin{equation}    \label{eq4.12}
\|u\|_{1,s}\leq C_1\kappa(2-s)\|u\|_{1,0},
\end{equation}
with the same constant $C_1$ as in \rmref{eq4.6} and 
\rmref{eq4.8} that depends on the total charge of 
the nuclei, on the number of electrons, and on the 
eigenvalue $\lambda$, but is independent of $s$, 
and with
\begin{equation}    \label{eq4.13}
\kappa(2-s)=
\frac{2}{\sqrt{\pi}}\frac{\Gamma((1-s)/2)}{\Gamma((2-s)/2)}
\leq\frac{2}{1-s}.
\end{equation}
\end{thrm}

\begin{proof}
By Lemma~\ref{lm4.5} and Lemma~\ref{lm4.6}, we 
can choose $\Omega$ so large that the assigned 
operator $PT(\lambda)$ is contractive as an 
operator from $\FH^1$ to $\FH^1$ and from $\X_0$ 
to $\X_0$. As an eigenfunction $u\in\FH^1$ for 
the eigenvalue $\lambda<0$ is a solution of the 
equation (\ref{eq4.2}), its high-frequency part 
$v=Pu$ in $\FH^1$ solves the equation (\ref{eq4.9}). 
Since $PT(\lambda)$ is contractive as an operator 
from $\FH^1$ to $\FH^1$, $v=Pu$ is the only 
solution of this equation in $\FH^1$ with $u-Pu$ 
given. This means that the high-frequency part 
$Pu$ of the eigenfunction possesses the 
representation
\begin{displaymath}
Pu=\sum_{k=1}^\infty (PT(\lambda))^k(u-Pu)
\end{displaymath}
in form of a Neumann series. The function $u-Pu$ 
is square integrable and has a bounded support, 
the ball of radius $\Omega$ around the origin. 
Therefore, it lies in $\X_0$. Thus, the Neumann 
series converges to an element in $\X_0$. 
That is, its high frequency part $Pu$ and the 
eigenfunction $u=Pu+(u-Pu)$ itself lie in $\X_0$.
Since $u=T(\lambda)u$ by (\ref{eq4.2}) and 
(\ref{eq4.4}), $u$ lies in $\X_s$ for $s<1$ by 
Lemma~\ref{lm4.4} and its norms (\ref{eq1.3}) 
satisfy the estimate (\ref{eq4.12}). 
The inequality (\ref{eq4.13}) follows from 
(\ref{eq2.11}) and the subsequent remark. 
\end{proof}

Now we have come to the end and reached the goal of 
our efforts. Translated back to the representation 
in position space, the theorem states that the 
solutions of the electronic Schr\"odinger equation 
for eigenvalues below the ionization threshold
$\Sigma\leq 0$ lie in the Barron spaces $\B^s$, 
$s<1$. The weighted $L_1$-norm 
\begin{equation}    \label{eq4.14}
\|u\|_{1,s}=\int_{\mathbb{R}^{3N}}(1+\|\omega\|^2)^{s/2}
|u(\omega)|\domega
\end{equation}
of their Fourier transforms does not grow faster 
than $\sim 1/(1-s)$ as $s$ goes to one, as with 
the hydrogen ground state. The example of the 
hydrogen ground state shows that the limit $s=1$ 
cannot be reached in general. The norms are 
likely to degenerate if the associated eigenvalue
approaches the ionization threshold.


\bibliographystyle{plain}
\bibliography{references}

\begin{thebibliography}{10}

\bibitem{Chen-Lu-Lu-Zhou}
Z.~Chen, L.~Lu, Y.~Lu, and S.~Zhou.
\newblock A regularity theory for static {S}chr{\"o}dinger equations on
  $\mathbb{R}^d$ in spectral {B}arron spaces.
\newblock {\em SIAM J. Math. Anal.}, 55:557--570, 2023.

\bibitem{DeRyck-Mishra}
T.~{De Ryck} and S.~Mishra.
\newblock Numerical analysis of physics-informed neural networks and related
  models in physics-informed machine learning.
\newblock {\em Acta Numerica}, 33:633--713, 2024.

\bibitem{Dus-Ehrlacher}
M.~Dus and V.~Ehrlacher.
\newblock Two-layers neural networks for {S}chr{\"o}dinger eigenvalue problems.
\newblock https://arxiv.org/abs/2409.01640.

\bibitem{Fournais-Hoffmann-Ostenhof-Sorensen}
S.~Fournais, M.~Hoffmann-Ostenhof, T.~Hoffmann-Ostenhof, and T.~{{\O}stergard
  S{\o}rensen}.
\newblock Sharp regularity estimates for {C}oulombic many-electron wave
  functions.
\newblock {\em Commun. Math. Phys.}, 255:183--227, 2005.

\bibitem{Fournais-Hoffmann-Ostenhof-Sorensen_2}
S.~Fournais, M.~Hoffmann-Ostenhof, T.~Hoffmann-Ostenhof, and T.~{{\O}stergard
  S{\o}rensen}.
\newblock Analytic structure of many-body {C}oulombic wave functions.
\newblock {\em Commun. Math. Phys.}, 289:291--310, 2009.

\bibitem{Kato}
T.~Kato.
\newblock On the eigenfunctions of many-particle systems in quantum mechanics.
\newblock {\em Commun. Pure Appl. Math.}, 10:151--177,~1957.

\bibitem{Kreusler-Yserentant}
H.-C. Kreusler and H.~Yserentant.
\newblock The mixed regularity of electronic wave functions in fractional order
  and weighted {S}obolev spaces.
\newblock {\em Numer. Math.}, 121:781--802, 2012.

\bibitem{Meng}
L.~Meng.
\newblock On the mixed regularity of {N}-body {C}oulombic wavefunctions.
\newblock {\em ESAIM:M2AN}, 57:2257--2282, 2023.

\bibitem{Stein}
E.~M. Stein.
\newblock {\em Singular Integrals and Differentiability Properties of
  Functions}.
\newblock Princeton University Press, Princeton, 1970.

\bibitem{Xu}
J.~Xu.
\newblock The finite neuron method and convergence analysis.
\newblock {\em Commun. Comput. Phys.}, 28:1707--1745, 2020.

\bibitem{Yserentant_2004}
H.~Yserentant.
\newblock On the regularity of the electronic {Schr\"{o}\-dinger} equation in
  {H}ilbert spaces of mixed derivatives.
\newblock {\em Numer. Math.}, 98:731--759, 2004.

\bibitem{Yserentant_2010}
H.~Yserentant.
\newblock {\em Regularity and Approximability of Electronic Wave Functions},
  volume 2000 of {\em Lecture Notes in Mathematics}.
\newblock Springer, Heidelberg Dordrecht London New York, 2010.

\bibitem{Yserentant_2011}
H.~Yserentant.
\newblock The mixed regularity of electronic wave functions multiplied by
  explicit correlation factors.
\newblock {\em ESAIM: M2AN}, 45:803--824, 2011.

\end{thebibliography}


\newpage

\section*{Appendix. The best possible constant in Lemma~\ref{lm2.3}}

\newenvironment{unnumberedlemma}
{\vspace{1.75ex}\noindent\bf Lemma.\it }{}{\vspace{1.0ex}}

As already mentioned in Section~\ref{sec2}, 
the best possible, minimum and not further 
improvable constant in Lemma~\ref{lm2.3} is 
$\kappa(\vartheta)/2$, with $\kappa(\vartheta)$ 
as in (\ref{eq2.10}). The norm of $\K$ seen as 
an operator from  $L_1$ to $B_s$, $s<-1$, is 
therefore $\kappa(-s)/2$. This results from 
the following properties of the inner integral 
from the proof.

\begin{unnumberedlemma}         
For $\vartheta>1$ and all $\eta\neq 0$ in 
$\mathbb{R}^3$, the inequality 
\begin{equation*}    
\int_{\mathbb{R}^3}\frac{1}{\|\omega-\eta\|^2}\,
\frac{1}{(1+\|\omega\|^2)^{\vartheta/2}}\,\domega \,< 
\int_{\mathbb{R}^3}\frac{1}{\|\omega\|^2}\,
\frac{1}{(1+\|\omega\|^2)^{\vartheta/2}}\,\domega
\end{equation*}
holds. As a function of $\eta$, the integral 
is continuous at $\eta=0$ and takes its 
maximum value there.
\end{unnumberedlemma}

\def \dphi    {\,\diff{\varphi}}
\def \dtheta  {\,\diff{\theta}}

\begin{proof}
The integral is invariant to orthogonal transformations  
of $\eta$. This means that we can restrict ourselves 
to vectors $\eta$ pointing into the direction of the 
third coordinate axis. For $\eta=a e_3$, $a\geq 0$, 
one gets 
\begin{displaymath}
\int_{\mathbb{R}^3}\frac{1}{\|\omega-\eta\|^2}\,
\frac{1}{(1+\|\omega\|^2)^{\vartheta/2}}\,\domega \,=
\int_0^\infty\!\int_0^{2\pi}\!\!\int_{-\pi/2}^{\pi/2}\!
\,\frac{1}{r^2+a^2-2\,ar\sin\theta}
\,\frac{1}{(1+r^2)^{\vartheta/2}}\;r^2\cos\theta\,\dtheta\dphi\dr
\end{displaymath}
in spherical coordinates, which reduces to 
the double integral
\begin{displaymath}
\int_{\mathbb{R}^3}\frac{1}{\|\omega-\eta\|^2}\,
\frac{1}{(1+\|\omega\|^2)^{\vartheta/2}}\,\domega
\,=\, 2\pi \int_0^\infty\bigg(
\int_{-1}^1\frac{1}{r^2+a^2-2\,ar\,t}\dt
\bigg)\frac{r^2}{(1+r^2)^{\vartheta/2}}\,\dr.
\end{displaymath}
Introducing the at $t=1$ weakly singular 
function
\begin{displaymath}
K(t)=\pi t\,\big(\ln\big((1+t)^2\big)-\ln\big((1-t)^2\big)\big),
\end{displaymath}
for $a>0$, that is, $\eta\neq 0$, this 
finally leads to the representation
\begin{displaymath}
\int_{\mathbb{R}^3}\frac{1}{\|\omega-\eta\|^2}\,
\frac{1}{(1+\|\omega\|^2)^{\vartheta/2}}\,\domega 
\,= \int_0^\infty\!
K\bigg(\frac{r}{a}\bigg)\frac{1}{(1+r^2)^{\vartheta/2}}\,\dr.
\end{displaymath}
For $0\leq t<1$ and $t>1$, respectively, the 
function $K(t)$ possesses the series expansions
\begin{displaymath}
K(t)=\,4\pi\,\sum_{k=0}^\infty\frac{1}{2k+1}\,t^{2k+2},
\quad
K(t)=\,4\pi\,\sum_{k=0}^\infty\frac{1}{2k+1}\,\frac{1}{t^{2k}}.
\end{displaymath}
The function $K_0(t)=K(t)-4\pi$ is therefore 
integrable over the entire positive real 
axis. The integrals
\begin{displaymath}
\int_0^1 K_0(t)\dt\,=-\,2\pi,
\quad
\int_1^\infty K_0(t)\dt\,=\,2\pi
\end{displaymath}
cancel each other out. Therefore, the above 
integral splits into the sum of the three 
terms
\begin{displaymath}
4\pi\int_0^\infty\!\frac{1}{(1+r^2)^{\vartheta/2}}\,\dr,
\;\;
\int_0^1 K_0(t)\delta(t)\dt,
\;\;
\int_1^\infty\! K_0(t)\delta(t)\dt,
\end{displaymath}
where $\delta(t)$ is the strictly decreasing 
function
\begin{displaymath}
\delta(t)=\frac{a}{(1+a^2t^2)^{\vartheta/2}} -
\frac{a}{(1+a^2)^{\vartheta/2}}.
\end{displaymath}
Because of $|\delta(t)|\leq a$, the last two 
terms tend to zero as $a$ goes to zero and
the integral to
\begin{displaymath}
4\pi\int_0^\infty\!\frac{1}{(1+r^2)^{\vartheta/2}}\,\dr
\,=\int_{\mathbb{R}^3} \frac{1}{\|\omega\|^2}\,
\frac{1}{(1+\|\omega\|^2)^{\vartheta/2}}\,
\domega.
\end{displaymath}
This proves the continuity at $\eta=0$. 
The function $K_0(t)$ increases on the interval
$0\leq t<1$ strictly from a value less than 
zero to infinity. Hence, it possesses a zero 
$t_0<1$. On the whole interval from zero to one
\begin{displaymath}
K_0(t)(\delta(t)-\delta(t_0))\,\leq\,0
\end{displaymath}
holds. The second of the three terms can 
therefore be estimated as
\begin{displaymath}
\int_0^1 K_0(t)\delta(t)\dt\,\leq\, 
\delta(t_0)\!\int_0^1 K_0(t)\dt
\end{displaymath}
and takes a value
less than zero. As $K_0(t)>0$ and $\delta(t)<0$
for $t>1$, the same holds for the third.
In summary,
\begin{displaymath}
\int_0^\infty\!K\bigg(\frac{r}{a}\bigg)\frac{1}{(1+r^2)^{\vartheta/2}}\,\dr
\,<\, 4\pi\int_0^\infty\!\frac{1}{(1+r^2)^{\vartheta/2}}\,\dr.
\end{displaymath}
This is the claimed inequality in 
spherical coordinates.
\end{proof}

\medskip

The starting point of the proof of 
Lemma~\ref{lm2.3} was the estimate
\begin{equation*}
\int_{\mathbb{R}^3}
\frac{1}{(1+\|\omega\|^2)^{\vartheta/2}}\,
|\K f(\omega)|\domega \,\leq\,
\frac{1}{2\pi^2}\int_{\mathbb{R}^3}\bigg(\,
\int_{\mathbb{R}^3}\frac{1}{\|\omega-\eta\|^2}\,
\frac{1}{(1+\|\omega\|^2)^{\vartheta/2}}\,
\domega\bigg)|f(\eta)|\deta
\end{equation*}
for rapidly decreasing functions 
$f:\mathbb{R}^3\to\mathbb{C}$, where equality 
holds for real-valued, nonnegative $f$.
Because of
\begin{equation*}
\frac{1}{2\pi^2}
\int_{\mathbb{R}^3}\frac{1}{\|\omega\|^2}\,
\frac{1}{(1+\|\omega\|^2)^{\vartheta/2}}\,\domega
=\frac{\kappa(\vartheta)}{2},
\end{equation*}
therefore the constant $\kappa(\vartheta)$ 
in the estimate (\ref{eq2.9}) can be 
replaced by $\kappa(\vartheta)/2$. Since
\begin{equation*}    
\int_{\mathbb{R}^3}\frac{1}{\|\omega-\eta\|^2}\,
\frac{1}{(1+\|\omega\|^2)^{\vartheta/2}}\,\domega 
\,\geq\, (1-\varepsilon)
\int_{\mathbb{R}^3}\frac{1}{\|\omega\|^2}\,
\frac{1}{(1+\|\omega\|^2)^{\vartheta/2}}\,\domega
\end{equation*}
for all $\eta$ in a sufficiently small 
neighborhood of the origin, conversely
\begin{equation*}    
\int_{\mathbb{R}^3}\frac{1}{(1+\|\omega\|^2)^{\vartheta/2}}\,
|\K f(\omega)|\domega
\,\geq\, 
(1-\varepsilon)\,\frac{\kappa(\vartheta)}{2}
\int_{\mathbb{R}^3}|f(\omega)|\domega
\end{equation*}
for real-valued, nonnegative functions 
$f$ that vanish outside this neigborhood.


\end{document}